\documentclass{imsart}

\RequirePackage[OT1]{fontenc}
\RequirePackage{amsthm,amsmath,natbib}
\RequirePackage[colorlinks,citecolor=blue,urlcolor=blue]{hyperref}

\RequirePackage{amssymb,mathtools }
\setcitestyle{authoryear,open={(},close={)}}

\startlocaldefs
\numberwithin{equation}{section}
\theoremstyle{plain}

\theoremstyle{remark}

\newcommand{\mathbbm}[1]{\text{\usefont{U}{bbm}{m}{n}#1}} 
\endlocaldefs

\begin{document}

\begin{frontmatter}
\title{Gram Charlier and Edgeworth expansion for sample variance}
\runtitle{Gram Charlier and Edgeworth expansion for sample variance}

\begin{aug}
\author{\fnms{Eric} \snm{Benhamou}\thanksref{a,e1}\ead[label=e1,mark]{eric.benhamou@aisquare.com,eric.benhamou@dauphine.eu}}

\address[a]{A.I. SQUARE CONNECT, 35 Boulevard d'Inkermann 92200 Neuilly sur Seine, France and LAMSADE, Université Paris Dauphine, Place du Maréchal de Lattre de Tassigny,75016 Paris, France.\printead{e1}}

\runauthor{E. Benhamou}

\affiliation{A.I. SQUARE CONNECT and LAMSADE, Paris Dauphine University}

\end{aug}

\begin{abstract}
In this paper, we derive a valid Edgeworth expansions for the 
Bessel corrected empirical variance when data are generated
by a strongly mixing process whose distribution can be arbitrarily. 
The constraint of strongly mixing process makes the problem not easy. 
Indeed, even for a strongly mixing normal process, the distribution is unknown.
Here, we do not assume any other assumption than a sufficiently fast decrease of
the underlying distribution to make the Edgeworth expansion convergent.
This results can obviously apply to strongly mixing normal process and provide
an alternative to the work of \cite{Moschopoulos_1985} and \cite{Mathai_1982}.

\end{abstract}


\end{frontmatter}


\section{Introduction}
Let $X_1, \ldots, X_n$ be a random sample and define the sample variance statistic as:
\begin{equation}\label{definition}
\bar{X}_n =\frac{1}{n}\sum_{i=1}^{n}X_i, \quad s_n^2 = \frac{1}{n-1}\sum_{i=1}^{n}(X_i - \bar{X}_n)^2, \quad X_n = (X_1, \ldots, X_n)^T
\end{equation}
where $\bar X_n $ is the empirical mean, $s_n^2$ the Bessel corrected empirical variance also called sample variance, and $X_n$ the vector of the full history of this random sample. We are interested in the distribution of the sample variance under very weak conditions, namely that it admits a valid Edgeworth expansion. \\\\
It is insightful to notice that even with the additional constraint of a multi dimensional Gaussian distribution $\mathcal{N}(0, \Sigma)$\footnote{with $Sigma$ 
arbitrarily} for the underlying random vector $X_n$, the distribution of the sample variance is not known. In this particular setting, the sample variance is the squared norm of a multi dimensional Gaussian and can be seen as the linear combination of independent but not Homoscedastic variables. Standard theory states that the sample variance of a collection of independent and identically distributed normal variables follows a chi-squared distribution. But in this particular case, the different variables $X_i - \bar{X}_n$ are not independent and the result can not apply. Though they marginally have the same variance (if we conditioned by $\bar{X}_n$), they are correlated with each other. \\\\
Hence, in general, the distribution of the sample variance in the normal case is not a known distribution. There is however one special case, where it is a $\chi^2$ distribution.
Remember that when $X \sim \mathcal{N}(0, \Sigma)$, $X^T B X$ has a $\chi^2$ distribution of degree $d$ the rank of $B$ if and only if $B \Sigma B = B$ (see for instance page 413 remark following theorem 1 of \cite{Hogg_1978}). Notice that in our case, the $B$ matrix is $B=I_n - \frac 1 n \mathbbm{1}_n^T \mathbbm{1}_n$, where $I_n$ is the identity matrix of size $n$ and $\mathbbm{1}_n$ is the vector of size $n$ filled with 1 as we can write

\begin{eqnarray}
s_n^2 &= &\frac{1}{n-1}\sum_{i=1}^{n}(X_i - \bar{X}_n)^2 \\
& =& \frac{1}{n-1} X^T (I_n - \frac 1 n \mathbbm{1}_n^T \mathbbm{1}_n) X 
\end{eqnarray}

we obtain that it is a $\chi^2$ distribution if and only if 
$$ (I_n - \frac 1 n \mathbbm{1}_n^T \mathbbm{1}_n) \Sigma (I_n - \frac 1 n \mathbbm{1}_n^T \mathbbm{1}_n) = I_n - \frac 1 n \mathbbm{1}_n^T \mathbbm{1}_n $$
In other cases, the sample variance is a linear combination of Gamma distribution, and one has to rely on approximations as explained in \cite{Moschopoulos_1985} and \cite{Mathai_1982}.  \\\\
This simple example explains the interest in deriving an approximation of the distribution of the sample variance by means of Gram Charlier and Edgeworth expansion.

\section{Gram Charlier and Edgeworth expansion}
\subsection{Key concepts}
Gram–Charlier expansion\footnote{named in honor of the Danish mathematician, Jørgen Pedersen Gram and the Swedish astronomer, Carl Charlier}, and Edgeworth expansion\footnote{named in honor of the Anglo-Irish philosopher, Francis Ysidro Edgeworth}, are series that approximate a probability distribution in terms of its cumulants. The series are the same but, they differ in the ordering of their terms. Hence the truncated series are different, as well as the accuracy of truncating the series. The key idea in these two series is to expand the characteristic function in terms of the characteristic function of a known distribution with suitable properties, and to recover the concerned distribution through the inverse Fourier transform.

In our case, a natural candidate to expand around is the normal distribution as the central limit theorem and its different extensions to non independent and non identically distributed variable state that the resulting distribution is a normal distribution (or in the most general case to truncated symmetrical  and $\alpha$ stable distributions \footnote{please see the extension of CLT in \cite{Gnedenko_Kolmogorov_1954}}.

Let denote by $\hat{f}$ (respectively $\hat{\phi}$ the characteristic function of our distribution whose density function is $f$ (respectively $\phi$), and $\kappa_j$ its cumulants (respectively $\gamma_j$). Cumulants definition state 

\begin{equation}
\hat{f}(t)= \exp\left[\sum_{j=1}^\infty\kappa_j\frac{(it)^j}{j!}\right] \quad \quad \text{and} \quad \quad \hat{\phi}(t)=\exp\left[\sum_{j=1}^\infty\gamma_j\frac{(it)^j}{j!}\right],
\end{equation}

which gives the following formal identity:
\begin{equation}
\hat{f}(t)=\exp\left[\sum_{j=1}^\infty(\kappa_j-\gamma_j)\frac{(it)^j}{j!}\right]\hat{\phi}(t)\,.
\end{equation}

Using Fourier transform property that say $(it)^j \hat{\phi}(t)$ is the Fourier transform of $(-1)^j[D^j\phi](-x)$, where $D$ is the differential operator with respect to $x$, we get the formal expansion:

\begin{equation}
f(x) = \exp\left[\sum_{j=1}^\infty(\kappa_j - \gamma_j)\frac{(-D)^j}{j!}\right]\phi(x)\,
\end{equation}

If $\phi$ is chosen as the normal density $\phi(x) = \frac{1}{\sqrt{2\pi}\sigma}\exp\left[-\frac{(x-\mu)^2}{2\sigma^2}\right]$  with mean and variance as given by $f$, that is, mean $\mu = \kappa_1$ and variance $\sigma^2 = \kappa_2$, then the expansion becomes

\begin{equation}
f(x) = \exp\left[\sum_{j=3}^\infty\kappa_j\frac{(-D)^j}{j!}\right] \phi(x),
\end{equation}

since $ \gamma_j=0$ for all $j  > 2$, as higher cumulants of the normal distribution are 0. By expanding the exponential and collecting terms according to the order of the derivatives, we arrive at the Gram–Charlier A series. Such an expansion can be written compactly in terms of Bell polynomials as

\begin{equation}
\exp\left[\sum_{j=3}^\infty\kappa_j\frac{(-D)^j}{j!}\right] = \sum_{j=0}^\infty B_j(0,0,\kappa_3,\ldots,\kappa_j)\frac{(-D)^j}{j!}.
\end{equation}

Since the j-th derivative of the Gaussian function $\phi$ is given in terms of Hermite polynomial as

\begin{equation}
\phi^{(j)}(x) = \frac{(-1)^j}{\sigma^j} He_j \left( \frac{x-\mu}{\sigma} \right) \phi(x),
\end{equation}

this gives us the final expression of the Gram-Charlier A series as

\begin{equation}
 f(x) = \phi(x) \sum_{j=0}^\infty \frac{1}{j! \sigma^j} B_j(0,0,\kappa_3,\ldots,\kappa_j)  He_j \left( \frac{x-\mu}{\sigma} \right).
\end{equation}

If we include only the first two correction terms to the normal distribution, we obtain
\begin{equation}
 f(x) = \phi(x) \left[1+\frac{\kappa_3}{3!\sigma^3}He_3\left(\frac{x-\mu}{\sigma}\right)+\frac{\kappa_4}{4!\sigma^4}He_4\left(\frac{x-\mu}{\sigma}\right) + R_5 \right]\,,
\end{equation}

with $He_3(x)=x^3-3x$ and $He_4(x)=x^4 - 6x^2 + 3$ and 
$$R_5(x) = \sum_{j=5}^\infty \frac{1}{j! \sigma^j} B_j(0,0,\kappa_3,\ldots,\kappa_j)  He_j \left( \frac{x-\mu}{\sigma} \right).$$

If in the above expression, the cumulant are function of a parameter $\frac 1 n$, we can rearrange terms by power of $\frac 1 n$ and find the Edgeworth expansion.

\subsection{Cumulant for weak conditions}
In order to derive our Gram Charlier or Edgeworth expansion, we need to compute in full generality our different cumulants. Using similar techniques as in \cite{Benhamou_2018_SampleVariance}, we can get the various cumulants as follows:

The first two cumulants are easy and given by:
\begin{small}
\begin{flalign}
  \kappa_1 = &  (\acute{\mu }_{2})_1- (\acute{\mu }_{1}^2)_1 & \\
\kappa_2  = &  \frac{A^2_{1,0}}{n-1} +  \frac{A^2_{0,1}}{n} + \frac{A^2_{1,1} } {(n-1) n } + R_2 & 
\end{flalign}
\end{small}

with the different numerator terms given by:
\begin{small}
\begin{flalign}
A^2_{1,0} = & - 4 (\acute{\mu }_1^4)_2+8  ( \acute{\mu }_1^2 \acute{\mu }_2)_2-(\acute{\mu }_{2}^2)_2  &\\
A^2_{0,1} = & (\acute{\mu }_4)_2-4 (\acute{\mu }_1 \acute{\mu }_3)_2 & \\
A^2_{1,1} = & 6( \acute{\mu }_1^4)_2  - 12 (\acute{\mu}_1^2 \acute{\mu }_2 )_2 +3 (\acute{\mu }_2^2)_2  & \\
R_2 = & (\acute{\mu }_2^2)_2 - 2 (\acute{\mu}_1^2 \acute{\mu }_2)_2  +( \acute{\mu }_1^4)_2 -  (\kappa_1)^2 & 
\end{flalign}
\end{small}
with the natural symmetric moment estimators whose expressions are provided in appendix section \ref{moment_notation_1_2}. The term $R_2$ is the second order rest and is equal to zero if the sample is i.i.d. For general case, this term does not cancel out and should be taken into account. It can be rewritten as

\begin{small}
\begin{flalign}
R_2 = & (\acute{\mu }_2^2)_2 - (\acute{\mu }_{2})_1^2  + 2 ((\acute{\mu }_{2})_1 (\acute{\mu }_{1}^2)_1-  (\acute{\mu}_1^2 \acute{\mu }_2)_2 ) +( \acute{\mu }_1^4)_2 -(\acute{\mu }_{1}^2)_1^2
\end{flalign}
\end{small}

The third cumulant is more involved and given by:
\begin{small}
\begin{flalign} 
\kappa_3  =   & \frac{A^3_{2,0}}{(n-1)^2} + \frac{A^3_{1,1}}{(n-1) n} + \frac{A^3_{0,2}}{n^2}+ \frac{A^3_{2,1}} {(n-1)^2 n}  +  \frac{A^3_{1,2}} {(n-1)^2 n}  \nonumber \\ 
	 & +  \frac{A^3_{2,2}} {(n-1)^2 n^2}  + R_3 
\end{flalign}
\end{small}

where the different numerator terms are given by:
\begin{small}
\begin{flalign}
A^3_{2,0} = & -40 (\acute{\mu }_1^6)_3
	+120 (\acute{\mu }_1^4 \acute{\mu }_2 )_3
	-56 (\acute{\mu }_1^3 \acute{\mu }_3 )_3
	-78 (\acute{\mu }_1^2 \acute{\mu }_2^2 )_3 
	+48 ( \acute{\mu }_1 \acute{\mu }_2 \acute{\mu }_3)_3 \nonumber \\ &
	+2 (\acute{\mu }_2^3)_3
  	-6 (\acute{\mu }_3^2)_3 \\
A^3_{1,1} = & 18 (\acute{\mu }_1^2 \acute{\mu }_4)_3
		-3 (\acute{\mu }_2 \acute{\mu }_4)_3 \\
A^3_{0,2} = & (\acute{\mu }_6 )_3 -6 (\acute{\mu }_1 \acute{\mu }_5)_3 \\
A^3_{2,1} = & 136 (\acute{\mu }_1^6)_3
		-408 (\acute{\mu }_1^4 \acute{\mu }_2 )_3
		+160 (\acute{\mu}_1^3 \acute{\mu }_3 )_3
		+288 (\acute{\mu }_1^2 \acute{\mu }_2^2 )_3
		-144 (\acute{\mu }_1 \acute{\mu }_2 \acute{\mu }_3 )_3 \nonumber \\ &
		-24 (\acute{\mu }_2^3)_3
		+12 (\acute{\mu}_3^2)_3 \\
A^3_{1,2} = &  15 (\acute{\mu }_2 \acute{\mu }_4)_3
		-30 (\acute{\mu }_1^2 \acute{\mu }_4)_3 \\
A^3_{2,2} = & -120 (\acute{\mu }_1^6)_3
	+360 (\acute{\mu }_1^4 \acute{\mu }_2 )_3
	-120  (\acute{\mu }_1^3 \acute{\mu }_3)_3
	-270 (\acute{\mu }_1^2 \acute{\mu }_2^2 )_3
	+120 ( \acute{\mu }_1 \acute{\mu }_2 \acute{\mu }_3)_3 \nonumber \\ &
	+30 (\acute{\mu }_2^3)_3
	-10 (\acute{\mu }_3^2 )_3 \\
R_3 = &  -3 \mu'_2\mu'_1 +2  (\kappa_1)^3 
\end{flalign}
\end{small}

with the natural symmetric moment estimators given in appendix section \ref{moment_notation_3}. The fourth cumulant is even more involved and given by:

\begin{small}
\begin{align}
\kappa_4 = &  \frac{A^4_{3,0}}{(n-1)^3}  + \frac{A^4_{2,1}}{(n-1)^2 n} + \frac{A^4_{1,2}}{(n-1)n^2} + \frac{A^4_{0,3}} {n^3}  +  \frac{A^4_{3,1}} {(n-1)^3 n} +  \frac{A^4_{2,2}} {(n-1)^2 n^2}  \nonumber \\
&+\frac{A^4_{1,3}}{(n-1) n^3} + \frac{A^4_{3,2}} {(n-1)^3 n^2} + \frac{A^4_{2,3}} {(n-1)^2 n^3} + \frac{A^4_{3,3}}{(n-1)^3 n^3} + R_4
\end{align}

with  the different numerator terms are given by:
\begin{flalign}
A^4_{3,0} = & -672 (\acute{\mu }_1^8 )_4
	+2688 (\acute{\mu }_1^6 \acute{\mu }_2 )_4
	-1216 (\acute{\mu }_1^5 \acute{\mu }_3 )_4
	-3120 (\acute{\mu }_1^4 \acute{\mu }_2^2 )_4
	+400  (\acute{\mu }_1^4 \acute{\mu }_4 )_4
	+2240 (\acute{\mu }_1^3 \acute{\mu }_2 \acute{\mu }_3 )_4  \nonumber \\
	& +960 (\acute{\mu }_1^2 \acute{\mu }_2^3 )_4
	-384 (\acute{\mu }_3^2  \acute{\mu }_1^2  )_4
	-480 (\acute{\mu }_1^2 \acute{\mu }_2 \acute{\mu }_4 )_4
	-624 (\acute{\mu }_1 \acute{\mu }_2^2 \acute{\mu }_3 )_4
	+144 (\acute{\mu }_1 \acute{\mu }_3 \acute{\mu }_4 )_4
	-6 (\acute{\mu }_2^4)_4 \nonumber \\ 
	& +96 (\acute{\mu }_2 \acute{\mu }_3^2)_4
	-3 (\acute{\mu }_4^2)_4
	+12 (\acute{\mu }_2^2 \acute{\mu}_4)_4  \\
A^4_{2,1} = & -128 (\acute{\mu }_1^3 \acute{\mu }_5 )_4
	+96 (\acute{\mu }_1 \acute{\mu }_2 \acute{\mu }_5 )_4
	-24 (\acute{\mu }_3 \acute{\mu }_5)_4 \\
A^4_{1,2} = & 32 (\acute{\mu }_1^2 \acute{\mu }_6)_4
	-4 (\acute{\mu }_2 \acute{\mu}_6)_4 \\
A^4_{0,3} = & (\acute{\mu }_8)_4
	-8 (\acute{\mu }_1 \acute{\mu }_7)_4 \\
A^4_{3,1} = &  3792 (\acute{\mu }_1^8 )_4
	-15168 (\acute{\mu }_1^6 \acute{\mu }_2 )_4
	+6144 (\acute{\mu }_1^5 \acute{\mu }_3 )_4
	+18144(\acute{\mu }_1^4 \acute{\mu }_2^2 )_4
   -1920 (\acute{\mu }_1^4 \acute{\mu }_4 )_4
   -11520 (\acute{\mu }_1^3 \acute{\mu }_2 \acute{\mu }_3 )_4  \nonumber \\ 
	& -6336 (\acute{\mu }_2^3  \acute{\mu }_1^2  )_4
	+1680 (\acute{\mu }_1^2 \acute{\mu }_3^2 )_4
	+2520 (\acute{\mu }_1^2 \acute{\mu }_2 \acute{\mu }_4 )_4
	+3600 (\acute{\mu }_1 \acute{\mu }_2^2 \acute{\mu}_3 )_4
   -624 (\acute{\mu }_1 \acute{\mu }_3 \acute{\mu }_4 )_4
   +234 (\acute{\mu }_2^4)_4 \nonumber \\ 
	& -432 (\acute{\mu }_2 \acute{\mu }_3^2)_4
   +33 (\acute{\mu}_4^2)_4
   -252 (\acute{\mu }_2^2 \acute{\mu }_4)_4 \\
A^4_{2,2} = & 400 (\acute{\mu }_1^3 \acute{\mu }_5 )_4
	-336 (\acute{\mu }_1 \acute{\mu }_2 \acute{\mu }_5 )_4
	+48 (\acute{\mu }_3 \acute{\mu }_5)_4 \\
A^4_{1,3} = & 28 (\acute{\mu }_2 \acute{\mu }_6)_4
	-56 (\acute{\mu }_1^2 \acute{\mu }_6)_4 \\
A^4_{3,2} = & -7440 (\acute{\mu }_1^8 )_4
	+29760 (\acute{\mu }_1^6 \acute{\mu }_2 )_4
	-10880 (\acute{\mu }_1^5 \acute{\mu }_3 )_4
	-36480 (\acute{\mu }_1^4 \acute{\mu }_2^2 )_4
   +3104 (\acute{\mu }_1^4 \acute{\mu }_4 )_4
   +20992 (\acute{\mu }_1^3 \acute{\mu }_2 \acute{\mu }_3 )_4 \nonumber \\ 
	& +13824 (\acute{\mu }_1^2 \acute{\mu }_2^3 )_4 
	-2752 (\acute{\mu }_1^2 \acute{\mu }_3^2 )_4
	-4368 (\acute{\mu }_1^2 \acute{\mu }_2 \acute{\mu }_4 )_4
	-7248 (\acute{\mu }_1 \acute{\mu }_2^2 \acute{\mu }_3 )_4
	+976 (\acute{\mu }_1 \acute{\mu }_3 \acute{\mu }_4 )_4
   -738 (\acute{\mu }_2^4)_4  \nonumber \\ 
	& +800 (\acute{\mu }_2 \acute{\mu }_3^2)_4
	-57 (\acute{\mu }_4^2)_4
	+612 (\acute{\mu }_2^2 \acute{\mu }_4)_4 \\
A^4_{2,3} = & -336 (\acute{\mu }_1^3 \acute{\mu }_5 )_4
	+336 (\acute{\mu }_1 \acute{\mu}_2 \acute{\mu }_5 )_4
	-56 (\acute{\mu }_3 \acute{\mu }_5)_4 \\
A^4_{3,3} = & 5040 (\acute{\mu }_1^8 )_4
	-20160 (\acute{\mu }_1^6 \acute{\mu }_2 )_4
	+6720 (\acute{\mu }_1^5 \acute{\mu }_3 )_4
   	+25200 (\acute{\mu }_1^4 \acute{\mu }_2^2 )_4
   	-1680 ( \acute{\mu }_1^4 \acute{\mu }_4 )_4 
	-13440 (\acute{\mu }_1^3 \acute{\mu }_2 \acute{\mu }_3 )_4 \nonumber \\  &
	 -10080 (\acute{\mu }_1^2 \acute{\mu }_2^3 )_4 
	+1680 (\acute{\mu }_1^2 \acute{\mu }_3^2 )_4
	+2520 (\acute{\mu }_1^2 \acute{\mu }_2 \acute{\mu }_4 )_4
	+5040 (\acute{\mu }_1 \acute{\mu}_2^2 \acute{\mu }_3 )_4
	-560 (\acute{\mu }_1 \acute{\mu }_3 \acute{\mu }_4 )_4 \nonumber \\  &
	+630 (\acute{\mu }_2^4)_4
	-560 (\acute{\mu }_2 \acute{\mu }_3^2)_4 
	+35 (\acute{\mu }_4^2)_4
	-420 (\acute{\mu }_2^2 \acute{\mu }_4 )_4
\end{flalign}
\end{small}

with the natural symmetric moment estimators given in appendix section \ref{moment_notation_4}.

\begin{proof} see appendix section \ref{proof_cumulant}
\end{proof}

\section{Conclusion}
In this paper, we have derived the most general formula for the Gram Charlier and the resulting Edgeworth expansion for the sample variance under very weak conditions. Our formula does not assume that the underlying sample is independent neither identically distributed. This formula can therefore be applied to strong mixing processes like sample of an auto regressive process of order 1 (AR(1)). It extends in particular the work of \cite{Mikusheva_2015}
\clearpage

\appendix
\section{Notations}\label{moment_notation}
\subsection{Empirical Moments of order 1 and 2 Notation}\label{moment_notation_1_2}
We adopt the following notations
\begin{itemize}
\item moment of order 1:

\begin{small}
\begin{flalign}
(\acute{\mu }_2)_1 = &  \frac{\sum_{i=1}^n \mathbb{E} \left[ X_i^2 \right] }{n}&  \\
(\acute{\mu }_1^2)_1 = & \frac{  \sum_{i \neq j} \mathbb{E}\left[X_i X_j \right] }{n(n-1)} & 
\end{flalign}
\end{small}

\item moments of order 2:
\begin{small}
\begin{flalign}
(\acute{\mu}_4)_2 = & \frac{\sum_{i=1}^n \mathbb{E} \left[X_i^4 \right]}{n} &\\
(\acute{\mu}_1 \acute{\mu}_3 )_2 = & \frac{ \sum_{i \neq j } \mathbb{E} \left[X_i^3 X_j \right]}{n (n-1)} & \\
(\acute{\mu}_2 ^2)_2 = & \frac{ \sum_{i \neq j} \mathbb{E} \left[X_i^2 X_j^2 \right]}{n (n-1)} & \\
(\acute{\mu}_1^2 \acute{\mu}_2)_2 = & \frac{\sum_{i \neq j \neq k} \mathbb{E} \left[X_i^2 X_j X_k \right]}{n (n-1)(n-2) } & \\
(\acute{\mu}_1^4)_2 = &  \frac{\sum_{i \neq j \neq k \neq l}\mathbb{E} \left[ X_i X_j X_k  X_l \right] }{n (n-1)(n-2) (n-3)} &  
\end{flalign}
\end{small}
\end{itemize}

\subsection{Empirical Moments of order 3 Notation}\label{moment_notation_3}
\begin{small}
\begin{flalign}
(\acute{\mu }_1^6)_3 = & \frac{\sum_{i \neq j \neq k \neq l \neq m \neq n}\mathbb{E} \left[ X_i X_j X_k  X_l X_m X_n\right] }{n (n-1)(n-2) (n-3)(n-4)(n-5)} &  \\
(\acute{\mu }_1^4 \acute{\mu }_2 )_3 = &  \frac{\sum_{i \neq j \neq k \neq l \neq m }\mathbb{E} \left[ X_i^2 X_j X_k  X_l X_m \right] }{n (n-1)(n-2) (n-3)(n-4)} &  \\
(\acute{\mu }_1^3 \acute{\mu }_3 )_3 = & \frac{\sum_{i \neq j \neq k \neq l  }\mathbb{E} \left[ X_i^3 X_j X_k  X_l \right] }{n (n-1)(n-2) (n-3)} &  \\
(\acute{\mu }_1^2 \acute{\mu }_2^2 )_3 = &  \frac{\sum_{i \neq j \neq k \neq l  }\mathbb{E} \left[ X_i^2 X_j^2 X_k  X_l \right] }{n (n-1)(n-2) (n-3)} &  \\
(\acute{\mu }_1^2 \acute{\mu }_4)_3 = &  \frac{\sum_{i \neq j \neq k  }\mathbb{E} \left[ X_i^4 X_j X_k \right] }{n (n-1)(n-2) } &  \\
(\acute{\mu }_1 \acute{\mu }_2 \acute{\mu }_3)_3 =&  \frac{\sum_{i \neq j \neq k  }\mathbb{E} \left[ X_i^3 X_j^2 X_k \right] }{n (n-1)(n-2) } &  \\
(\acute{\mu }_1 \acute{\mu }_5)_3 = &  \frac{\sum_{i \neq j }\mathbb{E} \left[ X_i^5 X_j\right] }{n (n-1)} &  \\
(\acute{\mu }_2^3)_3 = &  \frac{\sum_{i \neq j \neq k  }\mathbb{E} \left[ X_i^2 X_j^2 X_k^2 \right] }{n (n-1)(n-2) } &   \\
(\acute{\mu }_2 \acute{\mu }_4)_3 = &  \frac{\sum_{i \neq j }\mathbb{E} \left[ X_i^4 X_j^2 \right] }{n (n-1)} &  \\
(\acute{\mu }_3^2)_3 = &  \frac{\sum_{i \neq j }\mathbb{E} \left[ X_i^3 X_j^3 \right] }{n (n-1)} &  \\
(\acute{\mu }_6 )_3 = & \frac{\sum_{i }\mathbb{E} \left[ X_i^6 \right] }{n } & 
\end{flalign}
\end{small}

\subsection{Empirical Moments of order 4 Notation}\label{moment_notation_4}
\begin{small}
\begin{flalign}
(\acute{\mu }_1^8 )_4  = & \frac{\sum_{i \neq j \neq k \neq l \neq m \neq n \neq o \neq p}\mathbb{E} \left[ X_i X_j X_k  X_l X_m X_n X_o X_p\right] }{n (n-1)(n-2) (n-3)(n-4)(n-5)(n-6)(n-7)} &  \\
(\acute{\mu }_1^6 \acute{\mu }_2 )_4= & \frac{\sum_{i \neq j \neq k \neq l \neq m \neq n \neq o }\mathbb{E} \left[ X_i^2 X_j X_k  X_l X_m X_n X_o \right] }{n (n-1)(n-2) (n-3)(n-4)(n-5)(n-6)} &  \\
(\acute{\mu }_1^5 \acute{\mu }_3 )_4 = & \frac{\sum_{i \neq j \neq k \neq l \neq m \neq n }\mathbb{E} \left[ X_i^3 X_j X_k  X_l X_m X_n \right] }{n (n-1)(n-2) (n-3)(n-4)(n-5)} &  \\
(\acute{\mu }_1^4 \acute{\mu }_2^2)_4 = & \frac{\sum_{i \neq j \neq k \neq l \neq m \neq n }\mathbb{E} \left[ X_i^2 X_j^2 X_k  X_l X_m X_n \right] }{n (n-1)(n-2) (n-3)(n-4)(n-5)} &  \\
(\acute{\mu }_1^4 \acute{\mu }_4)_4 = & \frac{\sum_{i \neq j \neq k \neq l \neq m }\mathbb{E} \left[ X_i^4 X_j X_k  X_l X_m  \right] }{n (n-1)(n-2) (n-3)(n-4)} &  \\
(\acute{\mu }_1^3 \acute{\mu }_5 )_4 = & \frac{\sum_{i \neq j \neq k \neq l }\mathbb{E} \left[ X_i^5 X_j X_k  X_l  \right] }{n (n-1)(n-2) (n-3)} &  \\
(\acute{\mu }_1^3 \acute{\mu }_2 \acute{\mu }_3 )_4 = & \frac{\sum_{i \neq j \neq k \neq l \neq m }\mathbb{E} \left[ X_i^3 X_j^2 X_k  X_l X_m  \right] }{n (n-1)(n-2) (n-3)(n-4)} &  \\
(\acute{\mu }_1^2 \acute{\mu }_2 \acute{\mu }_4 )_4 = & \frac{\sum_{i \neq j \neq k \neq l }\mathbb{E} \left[ X_i^4 X_j^2 X_k  X_l  \right] }{n (n-1)(n-2) (n-3)} &  \\
(\acute{\mu }_1^2 \acute{\mu }_2^3)_4 = & \frac{\sum_{i \neq j \neq k \neq l \neq m }\mathbb{E} \left[ X_i^2 X_j^2 X_k^2  X_l X_m  \right] }{n (n-1)(n-2) (n-3)(n-4)} &  \\
(\acute{\mu }_1^2 \acute{\mu }_3^2)_4 = & \frac{\sum_{i \neq j \neq k \neq l }\mathbb{E} \left[ X_i^3 X_j^3 X_k  X_l  \right] }{n (n-1)(n-2) (n-3)} &  \\
(\acute{\mu }_1^2 \acute{\mu }_6)_4 = & \frac{\sum_{i \neq j \neq k }\mathbb{E} \left[ X_i^6 X_j X_k  \right] }{n (n-1)(n-2)} &  \\
(\acute{\mu }_1 \acute{\mu }_2^2 \acute{\mu }_3 )_4 = & \frac{\sum_{i \neq j \neq k \neq l }\mathbb{E} \left[ X_i^3 X_j^2 X_k^2  X_l  \right] }{n (n-1)(n-2) (n-3)} &  \\
(\acute{\mu }_1 \acute{\mu }_2 \acute{\mu }_5 )_4 = & \frac{\sum_{i \neq j \neq k }\mathbb{E} \left[ X_i^5 X_j^2 X_k  \right] }{n (n-1)(n-2)} &  \\
(\acute{\mu }_1 \acute{\mu }_3 \acute{\mu }_4)_4 = & \frac{\sum_{i \neq j \neq k }\mathbb{E} \left[ X_i^4 X_j^3 X_k  \right] }{n (n-1)(n-2)} &  \\
(\acute{\mu }_1 \acute{\mu }_7)_4 = & \frac{\sum_{i \neq j }\mathbb{E} \left[ X_i^7 X_j \right] }{n (n-1)} &  \\
(\acute{\mu }_2^2 \acute{\mu}_4)_4 = & \frac{\sum_{i \neq j \neq k }\mathbb{E} \left[ X_i^4 X_j^2 X_k^2  \right] }{n (n-1)(n-2)} &  \\
(\acute{\mu }_2^4)_4 = & \frac{\sum_{i \neq j \neq k \neq l }\mathbb{E} \left[ X_i^2 X_j^2 X_k^2  X_l^2  \right] }{n (n-1)(n-2) (n-3)} &  \\
(\acute{\mu }_2 \acute{\mu}_6)_4 = & \frac{\sum_{i \neq j }\mathbb{E} \left[ X_i^6 X_j^2 \right] }{n (n-1)} &  \\
(\acute{\mu }_2 \acute{\mu }_3^2)_4 = & \frac{\sum_{i \neq j \neq k }\mathbb{E} \left[ X_i^3 X_j^3 X_k^2 \right] }{n (n-1)(n-2)} & 
\end{flalign}
\begin{flalign}
(\acute{\mu }_3 \acute{\mu }_5)_4 = & \frac{\sum_{i \neq j }\mathbb{E} \left[ X_i^5 X_j^3 \right] }{n (n-1)} &  \\
(\acute{\mu }_4^2)_4 = & \frac{\sum_{i \neq j }\mathbb{E} \left[ X_i^4 X_j^4 \right] }{n (n-1)} &  \\
(\acute{\mu }_8)_4 = & \frac{\sum_{i}\mathbb{E} \left[ X_i^8  \right] }{n} &  
\end{flalign}
\end{small}

\section{Cumulant computation}\label{proof_cumulant}
First of all, the first four cumulants, denoted by $\kappa_i$ for $ii=1,\ldots,4$ are obtained through standard relationships with respect to moments denoted by $\mu'_i$ as follows:
\begin{align}
\kappa_1 = {} & \mu'_1 \\
\kappa_2 = {} & \mu'_2-{\mu'_1}^2 \\
\kappa_3 = {} & \mu'_3 -3 \mu'_2\mu'_1 +2{\mu'_1}^3 \\
\kappa_4 = {} & \mu'_4 -4 \mu'_1 \mu'_3 -3{\mu'_2}^2 +12{\mu'_1}^2 \mu'_2 -6{\mu'_1}^4
\end{align}
Hence we are left with computing the first four moments of the sample variance. 
The first two moments of the sample variance are easy to compute and given for instance in \cite{Benhamou_2018_SampleVariance}.
For the cumulant of order 3 and 4, we first compute the different moments and then regroup the terms.

\subsection{Third Moment for $s_n^2$}
Let us do some routine algebraic computation. We have
\begin{eqnarray}
s_n^6 &= & \frac{1}{n^3(n-1)^3} \left( (n-1) \sum_{i=1}^n X_i^2 - \sum_{i \neq j} X_i X_j  \right)^3 \\
& =& \frac{1}{n^3(n-1)^3} \left( (n-1)^3 (\sum_{i=1}^n X_i^2)^3 - 3 (n-1)^2 (\sum_{i=1}^n X_i^2)^2 ( \sum_{k \neq l} X_k X_l ) \right. \nonumber  \\
& & \left. + 3 (n-1) (\sum_{i=1}^n X_i^2) ( \sum_{k \neq l} X_k X_l )^2 + ( \sum_{k \neq l} X_k X_l )^3 \right)
\end{eqnarray}

Let us expand. The first expansion $(\sum_{i=1}^n X_i^2)^3 $ is easy and immediate:
\begin{eqnarray}
(\sum_{i=1}^n X_i^2)^3 & = &\sum_{i=1}^n X_i^6 + 3 \sum_{i \neq j} X_i^4 X_j^2 + \sum_{i \neq j \neq k} X_i^2 X_j^2 X_k^2 
\end{eqnarray}

In the expansion of $(\sum_{i=1}^n X_i^2)^2 (\sum_{j \neq k} X_j X_k)$, the possibilities are:
\begin{itemize}
\item $\sum_{i \neq j} X_i^5 X_j$
\item $\sum_{i \neq j \neq k} X_i^4 X_j X_k$  
\item $\sum_{i \neq j } X_i^3 X_j^3$ 
\item $\sum_{i \neq j \neq k} X_i^3 X_j^2 X_k$  
\item $\sum_{i \neq j \neq k \neq l} X_i^2 X_j^2 X_k X_l$ 
\end{itemize}

In the expansion of $(\sum_{i=1}^n X_i^2) (\sum_{j \neq k} X_j X_k)^2$, the possibilities are:
\begin{itemize}
\item $\sum_{i \neq j} X_i^4 X_j^2$  
\item $\sum_{i \neq j \neq k} X_i^4 X_j X_k$
\item $\sum_{i \neq j \neq k} X_i^3 X_j^2 X_k$  
\item $\sum_{i \neq j \neq k \neq l} X_i^3 X_j X_k X_l$  
\item $\sum_{i \neq j \neq k} X_i^2 X_j^2 X_k^2$  
\item $\sum_{i \neq j \neq k \neq l} X_i^2 X_j^2 X_k X_l$ 
\item $\sum_{i \neq j \neq k \neq l \neq m} X_i^2 X_j X_k X_l X_m$ 
\end{itemize}

In the expansion of $(\sum_{i \neq j} X_i X_j)^3$, the possibilities are:
\begin{itemize}
\item $\sum_{i \neq j} X_i^3 X_j^3$
\item $\sum_{i \neq j \neq k} X_i^3 X_j^2 X_k$
\item $\sum_{i \neq j \neq k} X_i^2 X_j^2 X_k^2$
\item $\sum_{i \neq j \neq k \neq l} X_i^3 X_j X_k X_l$
\item $\sum_{i \neq j \neq k \neq l} X_i^2 X_j^2 X_k X_l$
\item $\sum_{i \neq j \neq k \neq l \neq m} X_i^2 X_j X_k X_l X_m$
\item $\sum_{i \neq j \neq k \neq l \neq m \neq n} X_i X_j X_k X_l X_m$
\end{itemize}

which leads to 
\begin{footnotesize}
\begin{align} 
\mathbb{E}[ s_n^8 ] =   &  -\frac{(n-5) (n-4) (n-3) (n-2) \acute{\mu }_1^6}{(n-1)^2 n^2} 
+\frac{3 (n-5) (n-4) (n-3) (n-2) \acute{\mu }_2 \acute{\mu }_1^4}{(n-1)^2 n^2} \nonumber \\ &
+\frac{4 (n-3) (n-2) (3 n-5) \acute{\mu }_3 \acute{\mu }_1^3}{(n-1)^2 n^2}
-\frac{3 (n-3) (n-2) \left(n^2-6 n+15\right) \acute{\mu }_2^2 \acute{\mu }_1^2}{(n-1)^2 n^2} \nonumber \\ &
-\frac{3 (n-5) (n-2) \acute{\mu }_4 \acute{\mu }_1^2}{(n-1) n^2} 
-\frac{12 (n-2) \left(n^2-4 n+5\right) \acute{\mu }_2 \acute{\mu }_3 \acute{\mu }_1}{(n-1)^2 n^2}
-\frac{6 \acute{\mu }_5 \acute{\mu }_1}{n^2} \nonumber \\ &
-\frac{2 \left(3 n^2-6 n+5\right) \acute{\mu }_3^2}{(n-1)^2 n^2} 
+\frac{3 \left(n^2-2 n+5\right) \acute{\mu }_2 \acute{\mu }_4}{(n-1) n^2}  
+\frac{\acute{\mu }_6}{n^2} \nonumber \\ &
+\frac{(n-2) \left(n^3-3 n^2+9 n-15\right) \acute{\mu }_2^3}{(n-1)^2 n^2} 
\end{align} 
\end{footnotesize}

Regrouping all the terms leads to
\begin{align} 
\mathbb{E}[ s_n^6 ] =   & M^3_{0,0} +\frac{M^3_{1,0}}{n-1} + \frac{M^3_{0,1}}{n} + \frac{M^3_{2,0}}{(n-1)^2} + \frac{M^3_{1,1}}{(n-1) n} + \frac{M^3_{0,2}}{n^2} \nonumber \\ &
 + \frac{M^3_{2,1}} {(n-1)^2 n}  +  \frac{M^3_{1,2}} {(n-1) n^2}  +  \frac{M^3_{2,2}} {(n-1)^2 n^2}   
\end{align}

with
\begin{footnotesize}
\begin{align} 
M^3_{2,0} = & 
-60 \acute{\mu } _1^6+180 \acute{\mu } _1^4 \acute{\mu } _2 -68 \acute{\mu } _1^3 \acute{\mu } _3 -129 \acute{\mu } _1^2 \acute{\mu } _2^2 +60 \acute{\mu } _2 \acute{\mu } _3 \acute{\mu } _1+13 \acute{\mu } _2^3-6 \acute{\mu } _3^2    \\
M^3_{0,0} = & -\acute{\mu }_1^6+3 \acute{\mu}_1^4 \acute{\mu }_2 -3 \acute{\mu }_1^2 \acute{\mu }_2^2 +\acute{\mu }_2^3    \\
M^3_{2,2} = & -120 \acute{\mu }_1^6+360 \acute{\mu }_1^4 \acute{\mu }_2 -120
   \acute{\mu }_1^3 \acute{\mu }_3 -270 \acute{\mu }_1^2 \acute{\mu }_2^2 +120 \acute{\mu }_1 \acute{\mu }_2 \acute{\mu }_3 +30 \acute{\mu }_2^3-10
   \acute{\mu }_3^2   \\
M^3_{1,0} = &12 \acute{\mu } _1^6-36  \acute{\mu } _1^4 \acute{\mu } _2+12 \acute{\mu } _1^3 \acute{\mu } _3 +27  \acute{\mu } _1^2 \acute{\mu } _2^2 -12 \acute{\mu } _1 \acute{\mu } _2 \acute{\mu } _3 -3 \acute{\mu } _2^3   \\
M^3_{2,1} = & 154 \acute{\mu }_1^6-462 \acute{\mu }_2
   \acute{\mu }_1^4+172 \acute{\mu }_1^3 \acute{\mu }_3 +333 \acute{\mu }_1^2 \acute{\mu }_2^2 -156 \acute{\mu }_1 \acute{\mu }_2 \acute{\mu }_3 -33
   \acute{\mu }_2^3+12 \acute{\mu }_3^2 \\
M^3_{1,0} = & -3 \acute{\mu }_1^2 \acute{\mu }_4 +3 \acute{\mu }_2 \acute{\mu }_4 \\
M^3_{1,2} = & 15 \acute{\mu }_2 \acute{\mu }_4-30 \acute{\mu }_1^2 \acute{\mu }_4 \\
M^3_{1,1} = & 21 \acute{\mu }_1^2 \acute{\mu }_4-6 \acute{\mu }_2 \acute{\mu }_4 \\
M^3_{0,2} = &  \acute{\mu }_6-6 \acute{\mu }_1 \acute{\mu }_5
\end{align}
\end{footnotesize}

Using previous results in the relationship between cumulant and moment leads to the result. \qed

\subsection{Fourth Moment for $s_n^2$}
Let us do some routine algebraic computation. We have
\begin{eqnarray}
s_n^8 &= & \frac{1}{n^4(n-1)^4} \left( (n-1) \sum_{i=1}^n X_i^2 - \sum_{i \neq j} X_i X_j  \right)^4 \\
& =& \frac{1}{n^4(n-1)^4} \left( (n-1)^4 (\sum_{i=1}^n X_i^2)^4 - 4 (n-1)^3 (\sum_{i=1}^n X_i^2)^3 ( \sum_{k \neq l} X_k X_l ) +  6 (n-1)^2 \right. \nonumber  \\
& & \left. (\sum_{i=1}^n X_i^2)^2   \sum_{k \neq l} X_k X_l )^2 +  4 (n-1) (\sum_{i=1}^n X_i^2) ( \sum_{k \neq l} X_k X_l )^3  + ( \sum_{k \neq l} X_k X_l )^4 \right)
\end{eqnarray}

Again, one needs to expand all the terms and look at all the various possibilities to demonstrate the following relationship where we have regrouped against each of the symmetric empirical moment estimator

\begin{footnotesize}
\begin{align} 
\mathbb{E}[ s_n^8 ] =   & \frac{(n-7) (n-6) (n-5) (n-4) (n-3) (n-2) \acute{\mu }_1^8}{(n-1)^3 n^3} 
-\frac{4 (n-7) (n-6) (n-5) (n-4) (n-3) (n-2) \acute{\mu }_2 \acute{\mu }_1^6}{(n-1)^3 n^3} \nonumber \\ &
-\frac{8 (n-5) (n-4) (n-3) (n-2) (3 n-7) \acute{\mu }_3 \acute{\mu }_1^5}{(n-1)^3 n^3} 
+\frac{6 (n-5) (n-4) (n-3) (n-2) \left(n^2-10 n+35\right) \acute{\mu }_2^2 \acute{\mu }_1^4}{(n-1)^3 n^3} \nonumber \\ &
+\frac{2 (n-4) (n-3) (n-2) \left(3 n^2-30 n+35\right) \acute{\mu }_4 \acute{\mu }_1^4}{(n-1)^3 n^3} 
+\frac{16 (n-4) (n-3) (n-2) \left(3 n^2-20 n+35\right) \acute{\mu }_2 \acute{\mu }_3 \acute{\mu }_1^3}{(n-1)^3 n^3} \nonumber \\ &
+\frac{8 (n-3) (n-2) (3 n-7) \acute{\mu }_5 \acute{\mu }_1^3}{(n-1)^2 n^3} 
-\frac{4 (n-4) (n-3) (n-2) \left(n^3-9 n^2+45 n-105\right) \acute{\mu }_2^3 \acute{\mu }_1^2}{(n-1)^3 n^3} \nonumber \\ &
+\frac{8 (n-3) (n-2) \left(9 n^2-30 n+35\right) \acute{\mu }_3^2 \acute{\mu }_1^2}{(n-1)^3 n^3} 
-\frac{12 (n-3) (n-2) \left(n^3-9 n^2+35 n-35\right) \acute{\mu }_2 \acute{\mu }_4 \acute{\mu }_1^2}{(n-1)^3 n^3} \nonumber \\ &
-\frac{4 (n-7) (n-2) \acute{\mu }_6 \acute{\mu }_1^2}{(n-1) n^3} 
-\frac{24 (n-3) (n-2) \left(n^3-7 n^2+25 n-35\right) \acute{\mu }_2^2 \acute{\mu }_3 \acute{\mu }_1}{(n-1)^3 n^3} \nonumber \\ &
-\frac{8 (n-2) \left(3 n^3-21 n^2+45 n-35\right) \acute{\mu }_3 \acute{\mu }_4 \acute{\mu }_1}{(n-1)^3 n^3} 
-\frac{24 (n-2) \left(n^2-4 n+7\right) \acute{\mu }_2 \acute{\mu }_5 \acute{\mu }_1}{(n-1)^2 n^3} \nonumber \\ &
-\frac{8 \acute{\mu }_7 \acute{\mu }_1}{n^3} 
+\frac{(n-3) (n-2) \left(n^4-4 n^3+18 n^2-60 n+105\right) \acute{\mu }_2^4}{(n-1)^3 n^3} \nonumber \\ &
-\frac{8 (n-2) \left(3 n^3-15 n^2+35 n-35\right) \acute{\mu }_2 \acute{\mu }_3^2}{(n-1)^3 n^3} 
+\frac{\left(3 n^4-12 n^3+42 n^2-60 n+35\right) \acute{\mu }_4^2}{(n-1)^3 n^3} \nonumber \\ &
+\frac{6 (n-2) \left(n^4-4 n^3+16 n^2-40 n+35\right) \acute{\mu }_2^2 \acute{\mu }_4}{(n-1)^3 n^3} 
-\frac{8 \left(3 n^2-6 n+7\right) \acute{\mu }_3 \acute{\mu }_5}{(n-1)^2 n^3} \nonumber \\ &
+\frac{4 \left(n^2-2 n+7\right) \acute{\mu }_2 \acute{\mu }_6}{(n-1) n^3} 
+\frac{\acute{\mu }_8}{n^3}
\end{align}
\end{footnotesize}

Once this is done, one needs to collect all the terms and can get the final result. \qed

\clearpage

\bibliography{mybib}

\end{document}